\documentstyle{amsppt}
\input amstex
\pageno=1
\NoBlackBoxes
\NoRunningHeads

\def\reg{\operatorname{reg}}

\topmatter
\title
Irrational asymptotic behaviour of  Castelnuovo-Mumford regularity
\endtitle
\author S. Dale Cutkosky\thanks  partially supported by NSF
\endthanks\\
{\rm Department of Mathematics, University of Missouri\\
Columbia, MO 65211, USA}\\ \\
\endauthor
\endtopmatter

\document

\heading 1. Introduction \endheading \medskip

Let $\Cal I$ be an ideal sheaf on $\Bbb P^r$.
The  regularity  $\reg(\Cal I)$ of $\Cal I$ is
defined to be the  least integer $t$ such that
$H^i(\Bbb P^r,\Cal I(t-i)) = 0$ for all $i \ge 1$ (c.f. Mumford, Lecture 14  [M]).

This concept can be defined purely algebraically. 
Let $A = k[X_1,\ldots,X_r]$ be a polynomial ring over an arbitrary field
$k$. Let $I$ be a homogeneous ideal in $A$.  The
 regularity $\reg(I)$ of $I$ is defined to be the
maximum degree $n$ for which there is an index $j$ such that $H_{\frak
m}^j(I)_{n-j} \neq 0$, where $H_{\frak m}^j(I)$ denotes the $j$th local
cohomology module of $I$ with respect to the maximal graded ideal $\frak
m$ of $A$ (c.f Eisenbud and Goto [EG], Bayer and Mumford [BM]).

If $\Cal I$ is the sheafification of a homogeneous ideal $I$ of $A$, 
then $\reg(\Cal I)=\reg(\widetilde{I})$, where $\widetilde{I}$
is the saturation of $I$.  The saturation of $I$ is defined by
$$
\widetilde{I} = \{f\in A\mid\text{ for each }0\le i\le n,\text{ there exists }n_i>0\text{ such that } x_i^{n_i}f\in I\}
$$

The functions $\reg (I^n)$ and $\reg(\Cal I^n)$ are quite interesting. Some references on this and related problems are
[BEL], [BPV], [Ch], [CTV], [GGP], [HT], [Hu2], [K],  [S], [SS], [ST], [T].

Let $I$ be any homogeneous ideal of $A$.
In Theorem 1.1 of [CHT] it is shown that
 $\reg (I^n)$ is a linear
polynomial for all $n$ large enough.

Because of the examples given in [CS] showing the failure of the Riemann Roch problem to have a good solution
in general, it is perhaps to be expected that  $\reg(\Cal I^n)$ will also not have polynomial like behaviour.
However, there are some significant differences in these problems, and there are special technical difficulties involved
 in constructing a geometric  example calculating the regularity of
 powers of an ideal sheaf in projective space.
 
In this paper we construct examples of  ideal sheaves $\Cal I$ 
 of  nonsingular curves in $\Bbb P^3$, showing bizarre behaviour of $\reg(\Cal I^n)$.

If we take $a=9$, $b=1$, $c=1$ in Theorem 10, and make use of the comments following Lemma 8,
we get
$$
\text{reg}(\Cal I_C^r) = [r\sqrt{2}]+1+\sigma(r)
$$
for $r>0$, where $[r\sqrt{2}]$ is the greatest integer in $r\sqrt{2}$,
$$
\sigma(r) = \left\{\matrix 0&\text{ if }r=q_{2n}\text{ for some }n\in\Bbb N\\
1&\text{ otherwise}
\endmatrix\right.
$$
$q_m$ is defined recursively by $q_0=1, q_1=2$ and $q_n=2q_{n-1}+q_{n-2}$.
Note that  $r$ such that $\sigma(r)=0$ are quite sparse,
as
$$
q_{2n}\ge 3^n.
$$
In this example,  
$$\lim {\reg(\Cal I^n) \over n}=\sqrt{2}.$$

Example 10 uses the method of [C], which is to find a rational line which intersects the boundary of the cone of
effective curves on a projective surface in irrational points.

In [CHT] we give  examples showing that  $\reg(\Cal I^n)$
is at least not eventually a linear polynomial, or  a linear polynomial with
periodic coefficients.
 In particular, using a counter-example
to Zariski's Riemann-Roch problem in positive characteristic
[CS] we  construct an example of a union points in $\Bbb P^2$
 in positive characteristic
such that $\reg(\Cal I^n)$ is not
eventually a  linear polynomial with periodic coefficients (Example 4.3 [CHT]).
In this example,
$$
\reg(\Cal I^{5n+1}) = 
\left\{\matrix 29n + 7 &\text{ if }n\text{ is not a power of }p\\
29n+8&\text{ if }n \text{ is a power of }p.
\endmatrix\right.
$$

In [CHT] we consider the limits
$$\lim {\reg(I^n) \over n}\tag 1$$
and 
$$\lim {\reg(\Cal I^n) \over n}.\tag 2$$
It follows from Theorem 1.1 of [CHT] that (1) always has a limit,
which must be a natural number.  
In the examples given in [CHT], (2)  has a limit which is
a rational number. In the examples of Theorem 10,
(2) is
an irrational number.   
$$\lim {\reg(\Cal I^n) \over n}\not\in \Bbb Q.$$

It is an interesting problem to try to find examples where the limit
(2) does not exist.

I would like to thank Keiji Oguiso for suggesting the intersection form $q = 4x^2-4y^2-4z^2$,
and Olivier Piltant and Qi Zhang for helpful discussions on the material of
this paper.

\heading 2. Preliminaries on K3 surfaces \endheading \medskip

A K3 surface is a nonsingular complex projective surface $F$ such that
 $H^1(F,\Cal O_F) = 0$ and the canonical divisor $K_F$ of $F$ is trivial.
 
For a divisor $D$ on a complex variety $V$ we will denote
$$
h^i(D) = \text{dim}_{\Bbb C}H^i(V,\Cal O_V(D)).\tag 3
$$
We will write $D\sim E$ if $D$ is linearly equivalent to $E$.

On a K3 surface $F$, the Riemann-Roch Theorem (c.f. 2.1.1 [SD])
is
$$
\chi(D) = \frac{1}{2}D^2+2.\tag 4
$$
If $D$ is ample, we have 
$$
h^i(D) = 0\text{ for }i>0\tag 5
$$
 by the Kodaira Vanishing
Theorem (In dimension 2 this follows from Theorem IV.8.6 [BPV] and Serre duality).
Let $\text{Num}(F)=\text{Pic}(F)/\equiv$, where $D\equiv 0$ if and only
if $(D\cdot C)=0$ for all curves $C$ on $F$. 
$\text{Num}(F)\cong \text{Pic}(F)$, since $F$ is a K3 (c.f. (2.3) [SD]),
and $\text{Pic}(F)\cong \Bbb Z^n$, where $n\le 20$ (c.f. Proposition VIII.3.3 [BPV]).
$\text{Num}(F)$ has an intersection form $q$
induced by intersection of divisors, $q(D) =(D^2)$. Let 
$$
\text{N}(F)=\text{Num}(F)\otimes_{\Bbb Z}\Bbb R.
$$
Let $NE(F)$ be the smallest convex cone in $\text{N}(X)$ containing all
effective 1-cycles, $\overline{NE}(F)$ be the closure of $NE(F)$
in the metric topology. We can identify $\text{Pic}(F)$ with the integral points in $\text{N}(F)$.
 By abuse of notation, we will sometimes identify a divisor with
its equivalence class in $\text{N}(F)$.

Nakai's criterion for ampleness (c.f. Theorem 5.1 [H]) is
$D$ is ample if and only if $(D^2)>0$ and $(D\cdot C)>0$ for all\
curves $C$ on $F$.

$[\gamma]$ will denote the greatest integer in a rational number $\gamma$.

\proclaim{Theorem 1}(Theorem 2.9 [Mo]) 
For $\rho\le 11$, every even lattice of signature $(1,\rho-1)$
occurs as the Picard group  of a smooth, projective K3 surface.
\endproclaim

\heading 3. Construction of a Quartic Surface \endheading\medskip 
By Theorem 1,
$$
q = 4x^2-4y^2-4z^2
$$ 
is the intersection form of a projective K3 surface $S$.

\proclaim{Lemma 2} Suppose that $C$ is an integral curve on $S$. Then $C^2\ge 0$.
\endproclaim
\demo{Proof}
Suppose that $(C^2)<0$. Then $h^0(C)=0$, and $(C^2)=-2$ (c.f. the note after (2.7.1) [SD]).
But $4\mid (C^2)$ by our choice of $q$, so $(C^2)\ge 0$.
\enddemo

\proclaim{Theorem 3} Suppose that $D$ is ample on $S$. Then

\item{1.} $h^i(D)=0$ for $i>0$.

\item{2.} $h^0(D)=\frac{1}{2}(D^2)+2$.

\item{3.} $\mid D\mid$ is base point free.

\item{4.} There exists a nonsingular curve $\Gamma$ on $S$ such that $\Gamma\sim D$.

\item{5.} Let $g = \frac{1}{2}(D^2)+1$ be the genus of a general curve $\Gamma\in\mid D\mid$. Let
$\phi:S\rightarrow \Bbb P^g$ be the morphism induced by $\mid D\mid $. Then either

\item\item{a)} $\phi$ is an isomorphism onto a surface of degree $2g-2$, or

\item\item{b)} $\phi$ is a 2-1 morphism onto a rational surface of degree $g-1$.
\endproclaim

\demo{Proof}
1. and 2. are (5) and (4). We will prove 3. By 2. we can assume that $D$ is effective. Suppose that $C$ is an
irreducible fixed component of $\mid D\mid$. Then $h^0(C) =0$ and $(C^2)=-2$ by the note after (2.7.1) [SD].
But this is not possible by Lemma 2. Thus $\mid D\mid$ has no fixed component, and is thus base point free by
Corollary 3.2 [SD]. 4. Follows from the fact that $\mid D\mid$ is base point free, $(D^2)>0$ and by Bertini's Theorem.

We will now prove 5. By (4.1) [SD] and Theorem 6.1 [SD] either b) holds or $\phi$ is a birational morphism
 onto a normal surface of degree $2g-2$. Since $S$ has no curves with negative intersection number, $\phi$ has no
exceptional curves. By Zariski's Main Theorem, $\phi$ is an isomorphism.
\enddemo

By Lemma 2 we have
$$
\overline{NE}(S) = \{(x,y,z)\in \Bbb R^3 | 4x^2-4y^2-4z^2\ge 0,\,\,x\ge 0\}.
$$
By Nakai's criterion, cited in Section 2, a divisor $D$ is ample if and only if $D$ is in the interior of
$\overline{NE}(S)$. 

Let $D$ be a  divisor on $S$ with class $(1,0,0)\in NE(S)$. $D$ is ample since $D$ is in the
interior of $\overline{NE}(S)$. $(D^2)=4$ implies $h^0(D) =4$. Thus there exists an effective divisor $H$ 
such that $H=(1,0,0)$.  

\proclaim{Remark 4} If A nontrivial divisor $D=(a,b,c)\in \Bbb Z^3$ is linearly equivalent to an effective divisor,
then we must have $a>0$, since otherwise $D^2=q(D)<0$.
Thus we cannot have $H\sim D_1+D_2$ where $D_1$ and $D_2$ are effective divisors, both non trivial.
\endproclaim
\proclaim{Theorem 5}
$\mid H\mid$ embeds $S$ as a quartic surface in $\Bbb P^3$.
\endproclaim
\demo{Proof}
Consider the rational map
$$
\Phi: S\rightarrow \Bbb P^3
$$
induced by $\mid H\mid$. By Theorem 3, either 
$\Phi$ is an isomorphism onto
a quartic surface in $\Bbb P^3$, or else $\Phi$ is a 2-1 morphism onto a quadric surface in $\Bbb P^3$.

If $\phi$ is 2-1, then $\Phi$ is 2-1 onto a quadric surface $Q$ in $\Bbb P^3$, which must have, after suitable choice of
homogeneous coordinates, one of the forms
$$
x_0x_1-x_2x_3=0
$$
or 
$$
x_0x_1-x_2^2=0,
$$  
and $\Phi^*(\Cal O(1))\cong\Cal O_S(H)$. In either case there exists a reducible member of
$\Gamma(Q,\Cal O_Q(1))$, so that $\mid H\mid$ contains a reducible member, a contradiction to Remark 4.
\enddemo 
\vskip .2truein

Choose   $(a,b,c)\in\Bbb Z^3$ such that $a>0$, $a^2-b^2-c^2>0$ and $\sqrt{b^2+c^2}\not\in \Bbb Q$.
Since $(a,b,c)$ is in the interior of $\overline{NE}(S)$ there exists an ample divisor $A$ such that 
$A=(a,b,c)$. By Theorem 3, there exists a nonsingular curve $C$ on $S$ such that
$C\sim A$.
Consider the line
$tH-C$, $-\infty<t<\infty$ in $N(S)$. This line intersects the quadric $q=0$ in 2 irrational points
$$
\lambda_2=a+\sqrt{b^2+c^2}\text{ and }
$$
$$
\lambda_1=a-\sqrt{b^2+c^2}.
$$ 

\proclaim{Remark 6} Suppose that $m,r\in\Bbb N$. Then

1. $mH-rC\in \overline{NE}(S)$ and is ample if $r\lambda_2<m$.

2. $mH-rC\not\in \overline{NE}(S)$ and $rC-mH\not\in \overline{NE}(S)$ if $r\lambda_1<m<r\lambda_2$.

3. $rC-mH\in \overline{NE}(S)$ and is ample if $m<r\lambda_1$.

\endproclaim

We can  choose $(a,b,c)$
so that 
$$
7<\lambda_1<\lambda_2\tag 6
$$ 
and
$$
\lambda_2-\lambda_1>2.\tag 7
$$
By (4), for all $m,r\in\Bbb N$,
$$
\chi(mH-rC) = \frac{1}{2}(mH-rC)^2+2.\tag 8
$$
\proclaim{Theorem 7} Suppose that $m,r\in\Bbb N$. Then
$$
\align
h^1(mH-rC) &= \left\{\matrix 0&\text{ if } r\lambda_2<m\\
                           -\frac{1}{2}(mH-rC)^2-2 &\text{ if }r\lambda_1<m<r\lambda_2
             \endmatrix\right.\\
h^2(mH-rC) &= 0 \text{ if }r\lambda_1<m
\endalign
$$
\endproclaim
\demo{Proof}
By Remark 6, and (5), we have that $h^i(mH-rC)=0$ if $r\lambda_2<m$
and $i>0$. 

For $r\lambda_1<m<r\lambda_2$ we have by Remark 6 that both $mH-rC\not\in \overline{NE(S)}$ and $rC-mH\not\in \overline{NE(S)}$.
Thus 
$h^0(mH-rC) = 0$ and (by Serre duality) $h^2(mH-rC) = h^0(rC-mH)=0$. By (8),
$$
h^1(mH-rC) = -\chi(mH-rC) = -\frac{1}{2}(mH-rC)^2-2.
$$
If $m<r\lambda_1$,
$h^2(mH-rC) =h^0(rC-mH)=0$, by Remark 6.
 \enddemo

\proclaim{Lemma 8} For all $r\in\Bbb Z_+$, we must have either
$$
\align
h^1([r\lambda_2]H-rC)>0&\text{ or}\\
h^1([r\lambda_2]H-rC)=0&\text{ and }h^1(([r\lambda_2]-1)H-rC)\ne 0
\endalign
$$
\endproclaim
\demo{Proof} (7) implies that $r\lambda_1<[r\lambda_2]-1<[r\lambda_2]<r\lambda_2$.
Suppose that 
$$
h^1([r\lambda_2]H-rC)=h^1(([r\lambda_2]-1)H-rC)=0.
$$
Then by Theorem 7,
$$
-\frac{1}{2}([r\lambda_2]H-rC)^2-2=-\frac{1}{2}(([r\lambda_2]-1)H-rC)^2-2
$$
so that
$$
([r\lambda_2]H-rC)^2=(([r\lambda_2]-1)H-rC)^2
$$
and thus
$$
q([r\lambda_2]-ra,-rb,-rc)=q([r\lambda_2]-ra-1,-rb,-rc).
$$
We then have
$$
([r\lambda_2]-ra)^2=([r\lambda_2]-ra-1)^2
$$
so that $-2([r\lambda_2]-ra)+1 = 0$, a contradiction since $[r\lambda_2]-ra\in\Bbb Z$.
\enddemo
\vskip .2truein
The dependence of the vanishing  $h^1([r\lambda_2]H-rC)=0$ on $r$ is surprisingly subtle.
Set $d=b^2+c^2$. 
$$
\align
&h^1([r\lambda_2]H-rC)=0\text{ iff }\tag 9\\
&([r\lambda_2]-ra)^2-r^2b^2-r^2c^2=-1\text{ iff}\\
&[r\lambda_2]=ra+\sqrt{r^2(b^2+c^2)-1}\text{ iff}\\
&[r\sqrt{b^2+c^2}]=\sqrt{r^2(b^2+c^2)-1}\text{ iff}
\endalign
$$
There exists $s\in \Bbb N$ such that $r^2(b^2+c^2)-1=s^2$ iff
\vskip .2truein
There exists $s\in \Bbb N$ such that $s+r\sqrt d$ is a unit in $\Bbb Q(\sqrt d)$ with norm
$$
N(s+r\sqrt d)=-1.\tag 10
$$
\vskip .2truein

If $b=c=1$ so that $d=2$, we have that the integral solutions to (10) are
$$
s+r\sqrt d = \pm(1+\sqrt 2)^{2n+1}\text{ and } s+r\sqrt d = \pm(-1+\sqrt2)^{2n+1}
$$
for $n\in\Bbb N$ (c.f. Theorem 244 [HW]). As explained in the proof of  Theorem 244 in [HW], if $\frac{p_n}{q_n}$
are the convergents of
$$
\sqrt 2= 1+{1\over\displaystyle 2+
   {\strut 1\over\displaystyle 2 + \cdots}},
$$
then the solutions to (10) are $\pm(p_{2n}+q_{2n}\sqrt 2)$ and  $\pm(p_{2n}-q_{2n}\sqrt 2)$
for $n\in \Bbb N$. When we impose the condition that $r,s\in\Bbb N$, we get that the solutions to (9) are
$$
q_{2n},\,\, n\in \Bbb N
$$
 when $d=2$. $q_m$ are defined recursively by
$$
q_0=1,\,\, q_1=2,\,\, q_m=2q_{m-1}+q_{m-2}.
$$

\heading 4. Construction of the Example \endheading\medskip

Regard $S$ as a subvariety of $\Bbb P^3$, as embedded by $\mid H\mid$. There exists a hyperplane $H'$ of
$\Bbb P^3$ such that $H'\cdot S=H$.
Now let $\pi:X\rightarrow \Bbb P^3$ be the blow up of $C$. Let $E=\pi^*(C)$ be the exceptional surface, $\overline H=\pi^*(H')$.
Let $\overline S$ be the strict transform of $S$. Then $\pi^*(S) = \overline S+E$, $S\cong \overline S$ and 
$\overline S\cdot E = C$. Let $\Cal I_C$ be the ideal sheaf of $C$ in $\Bbb P^3$. Since $C$ is smooth, for $m,r\in\Bbb N$,
$$
\align
&\pi_*\Cal O_X(m\overline H-rE)\cong \Cal I_C^r(m)\tag 11\\
&R^i\pi_*\Cal O_X(m\overline H-rE)=0\text{ for }i>0,
\endalign
$$
(c.f. Proposition 10.2 [Ma]).
From the short exact sequence
$$
0\rightarrow \Cal O_X(-\overline S)\rightarrow \Cal O_X\rightarrow \Cal O_{\overline S}\rightarrow 0
$$ 
and the isomorphism $\Cal O_X(-4\overline H+E)\cong \Cal O_X(-\overline S)$
we deduce for all $m,r$  short exact sequences
$$
0\rightarrow \Cal O_X(m\overline H-rE)\rightarrow \Cal O_X((m+4)\overline H-(r+1)E)\rightarrow 
\Cal O_{S}((m+4)H-(r+1)C)\rightarrow 0\tag 12
$$
We have by (11), 
$$
\align
h^i(\Cal O_X(m\overline H))&=h^i(\Cal O_{\Bbb P^3}(mH')) = 0\text{ for  }m\in \Bbb N,\text{  and }0<i<3,\\
h^3(\Cal O_X(m\overline H))&=h^3(\Cal O_{\Bbb P^3}(mH')) = h^0((-m-4)H')=0\text{ for } m\ge 0.
\endalign
$$
In particular, we have 
$$
h^i(\Cal O_X(m\overline H))=0\text{ for }i>0,\,\,m\ge 0.\tag  13
$$

\proclaim{Theorem 9} Suppose that $m,r\in\Bbb N$. Then
$$
\align
h^1(m\overline H-rE)&=\left\{\matrix 0&m> r\lambda_2\\
                                     h^1(mH-rC)&\text{ if }m=[r\lambda_2]\text{ or }m=[r\lambda_2]-1
                        \endmatrix\right.\\
h^2(\overline mH-rE)&=0\text{ if }m>\lambda_1r\\
h^3(m\overline H-rE)&=0 \text{ if }m>4r.
\endalign
$$
\endproclaim

\demo{Proof}From (12) we have
$$
0\rightarrow \Cal O_X((n+4(t-1))\overline H-(t-1)E)\rightarrow \Cal O_X((n+4t)\overline H-tE)
\rightarrow \Cal O_S((n+4t)H-tC)\rightarrow 0.\tag 14
$$
From Theorem 7, (13) and induction applied to the long exact cohomology sequence of (14) we deduce that 
$h^i((4t+n)\overline H-tE)=0$ 
if $t\ge 0$, $n\ge 0$, $t\lambda_2<4t+n$, $i>0$. Thus
$h^i(m\overline H-rE) = 0$
if $i>0$, $m\ge 4r$, $r\lambda_2<m$. Since $4<\lambda_2$ by (6), we have
$$
h^i(m\overline H-rE) = 0\text{ if }i>0,r\lambda_2<m.
$$
From Theorem 7, (13) and (14), we also deduce
$h^2(m\overline H-rE)=0$ for $r\lambda_1<m$, since $4<\lambda_1$ by (6), and
$h^3(m\overline H-rE)=0$ for $4r\le m$. 

We now compute $h^1([r\lambda_2]\overline H-rE)$. By (6), $\lambda_2-4>1>r\lambda_2-[r\lambda_2]$ implies
$[r\lambda_2]-4>(r-1)\lambda_2$. We have a short exact sequence
$$
0\rightarrow \Cal O_X([r\lambda_2]-4)\overline H-(r-1)E)\rightarrow \Cal O_X([r\lambda_2]\overline H-rE)\rightarrow
\Cal O_{ S}([r\lambda_2]H-rC)\rightarrow 0
$$
Since $h^1(([r\lambda_2]-4)\overline H-(r-1)E) = h^2(([r\lambda_2]-4)\overline H-(r-1)E)=0$, we have
$$
h^1([r\lambda_2]\overline H-rE) = h^1([r\lambda_2]H-rC).
$$
A similar calculation shows that
$$
h^1(([r\lambda_2]-1)\overline H-rE) = h^1(([r\lambda_2]-1)H-rC).
$$
since (6) implies that $[r\lambda_2]-5>(r-1)\lambda_2$.
\enddemo
\vskip .2truein
Define
$$
\sigma(r) =\left\{\matrix 0&\text{ if }h^1([r\lambda_2]H-rC)=0\\
1&\text{ if }h^1([r\lambda_2]H-rC)\ne 0.
\endmatrix\right.
$$
\proclaim{Theorem 10}
$$
\text{reg}(\Cal I_C^r) = [r\lambda_2]+1+\sigma(r)
$$
for $r> 0$. In particular
$$\lim {\reg(\Cal I^r) \over r}=\lambda_2\not\in \Bbb Q.$$
\endproclaim
\demo{Proof}
For $i,m,r\ge 0$ we have, by (11),
$$
h^i(\Cal I_C^r(m))=h^i(m\overline H-rE).
$$
By Theorem 9 and Lemma 8, 
$$
h^1(\Cal I_C^r(t-1))=h^1((t-1)\overline H-rE) = \left\{\matrix 0& t\ge [r\lambda_2]+1+\sigma(r)\\
\ne 0&t=[r\lambda_2]+\sigma(r)
\endmatrix\right.
$$
By (7),
$([r\lambda_2]+1)-2 = [r\lambda_2]-1\ge r\lambda_2-2>\lambda_1r$ and 
$$
h^2(\Cal I_C^r(t-2))=h^2((t-2)\overline H-rE) = 0\text{ if }t\ge [r\lambda_2]+1.
$$
By (6),
$([r\lambda_2]+1)-3=[r\lambda_2]-2\ge r\lambda_2-3>4r$ and 
$$
h^3(\Cal I_C^r(t-3))=h^3((t-3)\overline H-rE) = 0\text{ if }t\ge [r\lambda_2] +1.
$$
\enddemo
\vskip .2truein
As a consequence of Theorem 10, we can construct the example stated in the Introduction.

\heading References \endheading \medskip

[BM] D. Bayer and D. Mumford, What can be computed in algebraic geometry? In:
D. Eisenbud and L. Robbiano (eds.), Computational Algebraic Geometry and
Commutative Algebra, Proceedings, Cortona 1991, Cambridge University Press,
1993, 1-48. \smallskip

[BEL] A. Bertram, L. Ein, and R. Lazarsfeld, Vanishing theorems, a theorem of
Severi, and the equations defining projective varieties, J. Amer. Math. Soc. 4
(1991), 587-602. \smallskip

[BPV] W. Barth, C. Peters, A. Van de Ven, Compact Complex Surfaces, Springer, 1984.
[CTV] M. Catalisano, N.V. Trung, and G. Valla, A sharp bound for
the regularity
index of fat points in general position, Proc. Amer. Math. Soc. 118 (1993),
717-724. \smallskip
[Ch] K. Chandler, Regularity of the powers of an ideal, Comm. Algebra 25
(1997), 3773-3776. \smallskip

[C] S.D. Cutkosky, Zariski decomposition of divisors on algebraic varieties,
Duke Math. J, 53 (1986),1 49-156. \smallskip

[CS] S.D. Cutkosky and V. Srinivas, On a problem of Zariski on dimensions of
linear systems, Ann. of Math. 137 (1993), 531-559. \smallskip

[CHT] S.D. Cutkosky, J. Herzog and N.V. Trung, Asymptotic behaviour of the Castelnuovo-Mumford regularity,
to appear in Compositio Math. \smallskip

[EG] D. Eisenbud and S. Goto, Linear free resolutions and minimal
multiplicities, J. Algebra 88 (1984), 107-184. \smallskip

[GGP] A. Geramita, A. Gimigliano, and Y. Pitteloud, Graded Betti numbers of
some embedded rational $n$-folds, Math. Ann. 301 (1995), 363-380. \smallskip

[HW] G.H. Hardy and E.M. Wright, An Introduction to the Theory of Numbers,
Oxford Univ Press, 1938.

[H] R. Hartshorne, Ample Subvarieties of Algebraic Varieties, Lect. Notes in Math. 156, Springer, 1970.
\smallskip

[HT] L. T. Hoa and N. V. Trung, On the Castelnuovo regularity and the
arithmetic degree of monomial ideals, Math. Z. 229 (1998), 5-9-537. \smallskip

[Hu2] C. Huneke, Uniform bounds in noetherian rings, Invent. Math. 107 (1992),
203-223. \smallskip

[K] V. Kodiyalam, Asymptotic behaviour of Castelnuovo-Mumford Regularity, to appear in Proc. AMS.
\smallskip

[Ma] H. Matsumura, Geometric Structure of the Cohomology Rings in Abstract Algebraic Geometry,
Mem. Coll. Sci. Univ. Kyoto (a) 32 (1959), 33-84.
\smallskip

[Mo] D.R. Morrison, On K3 surfaces with large Picard Number, Invent. Math. 75 (1984), 105-121.
\smallskip

[M] D. Mumford, Lectures on Curves on an Algebraic Surface, Annals of Math. Studies 59, 
Princeton Univ. Press, 1966.\smallskip

[SD] B. Saint-Donat, Projective Models of K3 Surfaces, Amer. J. Math. 96 (1974), 602-639.\smallskip

[S] I. Swanson, Powers of ideals, primary decompositions, Artin-Rees lemma and
regularity, Math. Ann. 307 (1997), 299-313. \smallskip

[SS] K.~Smith and I.~Swanson, Linear bounds on growth of associated primes,
Comm. Algebra 25 (1997), 3071-3079. \smallskip

[ST] B. Sturmfels, Four Counterexamples in Combinatorial Algebraic Geometry,
to appear in J. Algebra.\smallskip

[T] N.V. Trung, The Castelnuovo regularity of the Rees algebra and the
associated graded ring, Trans. Amer. Math. Soc., 350 (1998), 519-537.

\enddocument